\documentclass[12pt,a4paper]{article}
\usepackage[latin1]{inputenc}
\usepackage{amsmath}
\usepackage{amsfonts}
\usepackage{amssymb}
\usepackage{amsthm}
\usepackage{psfrag}
\usepackage{graphicx}
\usepackage{enumitem}
\usepackage{indentfirst}
\usepackage{bbm}
\newcommand{\lel}{\left\langle}
\newcommand{\rir}{\right\rangle}
\newcommand{\diag}{\text{diag}}
\newtheorem{theorem}{Theorem}[section]
\newtheorem{definition}[theorem]{Definition}
\newtheorem{lemma}[theorem]{Lemma}

\newtheorem{prop}[theorem]{Proposition}

\newtheorem{ass}[theorem]{Assumption}

\title{On Anticipated backward stochastic differential equations with Markov chain noise}
\date{}

\author{Zhe Yang \thanks{Department of Mathematics and Statistics, University of Calgary, 2500 University Drive NW, Calgary, AB, T2N 1N4, Canada, yangzhezhe@gmail.com}\and Robert J. Elliott \thanks{Haskayne School of Business, University of Calgary, 2500 University Drive NW, Calgary, AB, T2N 1N4, Canada, School of Mathematical Sciences, University of Adelaide, SA 5005, Australia, relliott@ucalgary.ca}}

\begin{document}
\maketitle
\begin{abstract}
In 2013, Lu and Ren \cite {luren} considered anticipated backward stochastic differential equations driven by finite state, continuous time Markov chain noise and established the existence and uniqueness of the solutions of these equations and a scalar comparison theorem. In this paper, we provide an estimate for their solutions and study the duality between these equations and stochastic differential delayed equations with Markov chain noise. Finally we derive another comparison theorem for
these solutions depending only on the two drivers.
\end{abstract}
\section{Introduction}
\indent In 2009, a new kind of
backward stochastic differential equations (BSDEs), called
anticipated BSDEs, was introduced in Peng and Yang \cite{PengYang} as follows:
\begin{equation}
\left \{ \begin{array}{ll} \hskip.1mmY_t =
\xi_T+\int^T_tf(s,Y_s,Z_s,Y_{s+\delta(s)}, Z_{s+\zeta(s)}
)ds-\int^T_t Z_s dB_s,&t\in [ 0,T];\\
\hskip.2mmY_{t}= \xi_t,&t\in[T,T+K];\\
Z_t=\eta_t,&t\in [ T,T+K].\end{array} \right. \nonumber
\end{equation}
They were motivated as the duality of stochastic differential delayed equations (SDDEs for short).
Here, $B$ is Brownian motion, $\xi_.,\eta_.$ are called the terminal conditions and $f$ is called the driver. Peng and Yang \cite{PengYang} provided the existence and uniqueness for the solutions of anticipated BSDEs under similar Lipschitz conditions and gave corresponding comparison results. In 2011 Xu \cite{xu1} obtained a necessary and sufficient condition for the
comparison theorem of multidimensional anticipated BSDEs. Xu also discussed a general comparison theorem for one-dimensional anticipated BSDEs in \cite{xu2}. In 2013 Yang and Elliott \cite{yang} gave a converse comparison theorem for anticipated BSDEs
and related non-linear expectations.\\
\indent In 2011, {\O}ksendal, Sulem and Zhang \cite{ksendal} studied existence and uniqueness theorems for time-advanced BSDEs driven both by Brownian motion
and compensated Poisson random measures. Wu, Wang and Ren \cite{wuwang} extended results of Peng and Yang \cite{PengYang} for anticipated BSDEs to non-Lipschitz generators. In 2013, Yang and Elliott \cite{yang2} derived the existence of solutions to one-dimensional
anticipated BSDEs with continuous coefficients, and showed the existence and comparison results of the minimal solutions. Zong \cite{zong} discussed the existence and uniqueness of the solutions of anticipated BSDEs driven by the Teugels martingales and established the corresponding comparison theorem.\\
\indent In 2012, van der Hoek and Elliott \cite{RE1} introduced a
market model where uncertainties are modeled by a finite state
Markov chain, instead of Brownian motion or related jump diffusions. The
Markov chain has a semimartingale representation involving a vector martingale $M=\{M_t\in\mathbb{R}^N,~t\geq 0\}$. BSDEs in this
framework were introduced by Cohen and Elliott \cite{Sam1} as
$$ Y_t = \xi + \int_t^T f(s, Y_s, Z_s) ds -\int_t^T  Z'_sdM_s
,~~~~~t\in[0,T].
$$ Here $M=\{M_t\in\mathbb{R}^N,~t\geq 0\}$ is a martingale coming from the semimartingale representation of the continuous time
Markov chain. Cohen and Elliott \cite{Sam2}, \cite{Sam3} gave some comparison results for
multidimensional BSDEs in the Markov Chain model.\\
\indent In 2013, Lu and Ren \cite {luren} discussed anticipated BSDEs driven by finite state, continuous time Markov chains:
\begin{equation*}
\left \{ \begin{array}{ll}
-dY_t = f(t,Y_t,Z_t,Y_{t+\delta(t)}, Z_{t+\zeta(t)} )dt- Z_t' dM_t,&t\in [ 0,T]; \\
\hskip.57cmY_{t}= \xi_t,&t\in[T,T+K];\\
\hskip.54cmZ_t=\eta_t,&t\in [ T,T+K].\end{array} \right.
\end{equation*}
 In the same paper, they established the existence and uniqueness of the solutions to this kind of equation.\\
\indent In this paper, we provide more properties of solutions to anticipated BSDEs with Markov chain noise. First we study how to bound the solutions by the terminal conditions and the driver. Then we deduce there exists a duality between these equations and stochastic differential delayed equations (SDDEs) on Markov chains. This means anticipated BSDEs with Markov chain noise exist naturally.\\
\indent Lu and Ren \cite{luren} also established a comparison theorem for
one-dimensional anticipated BSDEs on Markov chains, based on the comparison result for BSDEs in Cohen and Elliott \cite{Sam3}.
They used conditions involving not only the two drivers but also
the two solutions. We shall provide a comparison result involving
conditions only on the two drivers. This means the comparison result is easier to apply. For example, the penalization of reflected anticipated BSDEs on Markov chains and the converse comparison theorem for anticipated BSDEs on Markov chains can be established using our comparison result.\\
\indent The paper is organized as follows. In Section 2, we introduce the model and give some preliminary results. Section 3 provides a new proof of the solutions to anticipated BSDEs on Markov chains and an estimate of the solutions. In Section 4 we show the duality between these equations and SDDEs on Markov chains. We establish in Section 5 a comparison result for one-dimensional anticipated BSDEs with Markov chain noise.
\section{The  Model and Some Preliminary Results}\label{prelim}
 Consider a finite state Markov chain. Following the papers
\cite{RE1} and \cite{RE2} of van der Hoek and Elliott, we assume the
finite state Markov chain $X=\{X_t, t\geq 0 \}$ is defined on the
probability space $(\Omega,\mathcal{F},P)$ and the state space of
$X$ is identified with the set of unit vectors $\{e_1,e_2\cdots,e_N\}$ in
$\mathbb{R}^N$, where $e_i=(0,\cdots,1\cdots,0) ' $ with 1 in the
$i$-th position. Take $\mathcal{F}_t=\sigma\{X_s ; 0\leq s \leq t\}$ to be
the $\sigma$-algebra  generated by the Markov process $X=\{X_t\}$
and $\{\mathcal{F}_t\}$ to be its completed natural filtration. Since $X$ is a right continuous with left
limits (written RCLL) jump-process, then the filtration $\{\mathcal{F}_t\}$ is also
right-continuous. The  Markov chain has the semimartingale
representation:
\begin{equation}\label{semimartingale}
X_t=X_0+\int_{0}^{t}A_sX_sds+M_t.
\end{equation}
Here, $A=\{A_t, t\geq 0 \}$ is the rate matrix of the chain $X$ and
$M$ is a vector martingale (See Elliott, Aggoun and Moore
\cite{RE4}).
We assume the elements $A_{ij}(t)$ of $A=\{A_t, t\geq 0 \}$ are bounded. Then the martingale $M$ is square integrable.\\
\indent Denote by $[X,X]$ the optional quadratic variation of $X$, which is a $N\times N$ matrix process and $\lel X, X\rir$, the unique predictable $N\times N$ matrix process such that
$[X,X]-\lel X,X \rir$ is a matrix valued martingale and write $L$ for the matrix martingale process where:
$$
L_t = [X,X]_t - \lel X,X\rir_t, \quad t \in [0,T].
$$
It is shown in \cite{Sam1} that:
\begin{equation}\label{3}
\lel X, X\rir_t =  \int_0^t \diag(A_sX_s) ds - \int_0^t \diag(X_s)
A'_s ds - \int_0^t A_s \diag(X_s) ds.
\end{equation}
\indent For $n\in\mathbb{N}$, denote for $\phi \in \mathbb{R}^n$, the Euclidean norm $|\phi|_{n}=\sqrt{\phi'\phi}$ and for $\psi \in \mathbb{R}^{n\times n}$, the matrix norm $\|\psi\|_{n\times n}=\sqrt{Tr(\psi'\psi)}$.\\
Let $\Psi$ be the matrix
\begin{equation}\label{Psi}\Psi_t = \diag(A_tX_{t-})- \diag(X_{t-})A'_t - A_t \diag(X_{t-}).
\end{equation}
Then $d\langle X,X\rangle_t=\Psi_tdt.$ For any $t>0$, Cohen and Elliott \cite{Sam1,Sam3}, define the semi-norm $\|.\|_{X_t}$, for
$C, D \in \mathbb{R}^{N\times K}$ as:
\begin{align*}
\lel C, D\rir_{X_t} & = Tr(C' \Psi_tD), \\[2mm]
\|C\|^2_{X_t} & = \lel C, C\rir_{X_t}.
\end{align*}
We only consider the case where $C \in \mathbb{R}^N$, hence we
introduce the semi-norm $\|.\|_{X_t}$ as:
\begin{align*}
\nonumber
\lel C, D\rir_{X_t} & = C' \Psi_t D, \\[2mm]
\|C\|^2_{X_t} & = \lel C, C\rir_{X_t}.
\end{align*}
It follows from equation \eqref{3} that
\[\int_t^T \|C\|^2_{X_s} ds = \int_t^T  C' d\lel X, X\rir_s C.\]
\indent The following lemma comes from Yang, Ramarimbahoaka and Elliott \cite{zhedim}.
\begin{lemma}\label{normbound}
For any $B\in \mathbb{R}^N$,
$$ ~~~~\|B\|_{X_t} \leq \sqrt{3m} |B|_N, ~~\text{ for any }t\in[0,T],$$
where $m>0$ is the bound of $\|A_t\|_{N\times N}$, for any $t\in[0,T]$.
\end{lemma}
\indent Lemma \ref{Z2} is Lemma 3.1 in Cohen and Elliott \cite{Sam3}.
\begin{lemma}\label{Z2}
For $Z$, a predictable process in $\mathbb{R}^N$, verifying:
 \[E \left[ \int_0^t \|Z_u\|^2_{X_u} du\right] < \infty,\]
we have:
\begin{equation*}
 E \left[\left(\int_0^t  Z'_{u} dM_u  \right)^2\right] = E \left[ \int_0^t \|Z_u\|^2_{X_u} du\right].
\end{equation*}
\end{lemma}
 \indent Denote by $\mathcal{P}$, the $\sigma$-field generated by the processes defined on $(\Omega, P, \mathcal{F})$ which are predictable with respect to the filtration $\{\mathcal{F}_t\}_{t \in [0,\infty)}$. For any $t,s,r\in[0,\infty)$, $t\leq r \leq s$, consider the following spaces: \\[2mm]
$ L^2(\mathcal{F}_t;\mathbb{R}): =\{ \xi:~\xi$ is an $\mathbb{R} \text{-valued}~ \mathcal{F}_t $-measurable random variable such that \\
$~~~~~~~~~~~~~~~~~~~~~~E[|\xi|^2]< +\infty\};$\\[2mm]
$L^2_{\mathcal{F}}(t,s;\mathbb{R}): =\{\phi_.:[t,s]\times\Omega\rightarrow\mathbb{R};~ \phi_.~\text{is an adapted and RCLL process with}\\~~~~~~~~~~~~~~~~~~~~E[\int^s_t|\phi(t)|^2dt]<+\infty\}$;\\[2mm]
$H^2(t,s;\mathbb{R}^N)=\{\phi_.:[t,s]\times\Omega\rightarrow\mathbb{R}^N;~ \phi_.\in\mathcal{P}~\mbox{with  }E[\int^s_t\|\phi(t)\|_{X_t}^2dt]<+\infty\};$\\[2mm]
$ H^2(\mathcal{F}_r;\mathbb{R}^N): =\{ \varphi_r $ is an $\mathbb{R}^N \text{-valued}~ \mathcal{F}_r \text{-measurable random variable with~} \varphi_.$\\$~~~~~~~~~~~~~~~~~~~~~\in H^2(t,s;\mathbb{R}^N)\}.$\\[2mm]
\indent Consider the following one-dimensional BSDE with the Markov chain noise:
\begin{equation}\label{BSDEMC}
Y_t = \xi + \int_t^T f(s, Y_s, Z_s ) ds -\int_t^T  Z'_{s} dM_s
,~~~~~t\in[0,T].
\end{equation}
Here the terminal condition $\xi$ and the coefficient $f$ are known. Lemma \ref{existence} (Theorem 6.2 in Cohen and Elliott \cite{Sam1})
gives the existence and uniqueness result for solutions to the BSDEs
driven by Markov chains:
\begin{lemma}\label{existence}
Assume $\xi \in L^2(\mathcal{F}_T)$, the
function $f: \Omega \times [0, T] \times \mathbb{R} \times
\mathbb{R}^N \rightarrow \mathbb{R}$ satisfies a Lipschitz
condition, in the sense that there exists two constants $l_1, l_2>0$  such
that $P$-a.s. for each $y_1,y_2 \in \mathbb{R}$ and $z_1,z_2 \in
\mathbb{R}^{N}$, $t\in[0,T]$,
\begin{equation}\label{Lipchl}
|f(t,y_1,z_1) - f(t, y_2, z_2)| \leq l_1 |y_1-y_2| + l_2\|z_1
-z_2\|_{X_t},
\end{equation}
and for each $(y,z)\in \mathbb{R} \times \mathbb{R}^N$, the process $(f(t,y,z))_{t\in [0,T]}$ is predictable. We also assume $f$ satisfies
\begin{equation}\label{finite}
 E [ \int_0^T |f(t,0,0)|^2 dt] <\infty.
\end{equation}
Then there exists a solution $(Y, Z)\in L^2_{\mathcal{F}}(0,T;\mathbb{R})\times P^2_{\mathcal{F}}(0,T;\mathbb{R}^N)$
to BSDE (\ref{BSDEMC}). Moreover, this solution is
unique among $(Y,Z) \in L^2_{\mathcal{F}}(0,T;\mathbb{R})\times P^2_{\mathcal{F}}(0,T;\mathbb{R}^N)$ and up to indistinguishability for $Y$ and equality $d\langle
X,X\rangle_t$ $\times\mathbb{P}$-a.s. for $Z$.
\end{lemma}
\indent Campbell and Meyer \cite{campbell} gave the following definition:
\begin{definition}\label{defMoore}
The Moore-Penrose pseudoinverse of a square matrix $Q$ is the matrix $Q^{\dagger}$ satisfying the properties:\\[2mm]
  1) $QQ^{\dagger}Q = Q$ \\[2mm]
  2) $Q^{\dagger}QQ^{\dagger} = Q^{\dagger}$ \\[2mm]
  3) $(QQ^{\dagger})' = QQ^{\dagger}$ \\[2mm]
  4) $(Q^{\dagger}Q)'=Q^{\dagger}Q.$
\end{definition}
\indent Recall the matrix $\Psi_.$ given by \eqref{Psi}. The following lemma is Lemma 3.3 in Cohen and Elliott \cite{Sam3}.
\begin{lemma}\label{Psibound} For all $t$, both $\Psi$ and $\Psi^{\dagger}$ are bounded.
\end{lemma}
 \indent We adapt Lemma 3.5 in Cohen and Elliott \cite{Sam3} for our framework as follows:
\begin{lemma}\label{sam35}
For any driver satisfying \eqref{Lipchl} and \eqref{finite}, for any $Y_.$ and $Z_.$,
\[P(f(t,Y_{t-},Z_t) = f(t,Y_{t-}, \Psi_t\Psi_t^{\dagger} Z_t), ~\text{ for all}~ t \in [0,+\infty])=1\]
and
\[\int_0^t Z'_s dM_s = \int_0^t (\Psi_s\Psi^{\dagger}_s Z_s)' dM_s.\]
Therefore, without any loss of generality, we shall assume $Z_.=(\Psi\Psi^{\dagger} Z)_.$.
\end{lemma}
\begin{ass}\label{ass111}
Assume the Lipschitz constant $l_2$ of the driver $f$ given in \eqref{Lipchl} satisfies  $$~~~~~l_2\|\Psi_t^{\dagger}\|_{N \times N} \sqrt{6m}<1, ~~~\text{ for any }~t \in [0,T],$$ where $\Psi_.$ is given in \eqref{Psi} and $m>0$ is the bound of $\|A_t\|_{N\times N}$, for any $t\in[0,T]$.
\end{ass}
\indent The following lemma, which is a comparison result for BSDEs driven by a Markov chain, is found in Yang, Ramarimbahoaka and Elliott \cite{zhedim}.
\begin{lemma} \label{CT'} For $i=1,2,$ suppose $(Y_.^{(i)},Z_.^{(i)})$ is the solution of the
BSDE:
$$Y^{(i)}_t = \xi_i + \int_t^T f_i(s, Y^{(i)}_s, Z^{(i)}_s ) ds
- \int_t^T (Z_{s}^{(i)})' dM_s,\hskip.4cmt\in[0,T].$$
Assume $\xi_1,\xi_2\in L^2(\mathcal{F}_T;\mathbb{R})$, and $f_1,f_2:\Omega \times [0,T]\times \mathbb{R}\times \mathbb{R}^N \rightarrow \mathbb{R}$ satisfy conditions such that the above two BSDEs have unique solutions. Moreover assume $f_1$ satisfies \eqref{Lipchl} and Assumption \ref{ass111}.
If $\xi_1 \leq \xi_2 $, a.s. and $f_1(t,Y_t^{(2)}, Z_t^{(2)}) \leq f_2(t,Y_t^{(2)}, Z_t^{(2)})$, a.e., a.s., then
$$P( Y_t^{(1)}\leq Y_t^{(2)},~~\text{ for any } t \in [0,T])=1.$$
\end{lemma}
\indent The following lemma which gives the duality between the solutions to linear BSDEs and linear SDEs is Theorem 2 in \cite{Sam2}, adapted for our one-dimensional case with Markov chain noise:
\begin{lemma}(Linear BSDEs)\label{Duality}
Let $(\eta, \mu)$ be a $du \times P-a.s.$ bounded $(\mathbb{R}^{1\times N}, \mathbb{R})$ valued predictable process, $g \in P^2_{\mathcal{F}}(0,T,\mathbb{R})$ and $\xi \in L^2(\mathcal{F}_T)$. Then the linear BSDE given by
\begin{align*}
Y_t = \xi + \int_t^T (\mu_sY_s+\eta_sZ_s+g_s) ds - \int_t^T Z'_s dM_s, ~~~ t \in [0,T]
\end{align*}
has a unique solution $(Y,Z) \in L^2_{\mathcal{F}}(0,T;\mathbb{R})\times P^2_{\mathcal{F}}(0,T;\mathbb{R}^N) $, (up to appropriate sets of measure zero). Furthermore, if for all $s \in [t,T]$
\begin{equation}\label{sam10}
1 + \eta_s \Psi_s^{\dagger}(e_j-X_{s-})
\end{equation} is non-zero (invertible for the multi-dimensional case) for all $j$ such that $e'_j A_sX_{s-}>0$, except possibly on some evanescent set, then $Y$ is given by the explicit formula
\begin{equation}\label{sam11}
Y_t = E [ \xi U_T + \int_t^T g_s U_s ds | \mathcal{F}_t]
\end{equation}
up to indistinguishability. Here $U$ is the solution to the one-dimensional SDE:
$$
\begin{cases}
dU_s= U_s \mu_s ds + U_{s-} \eta_{s} (\Psi_s^{\dagger})'dM_s, ~~ s \in [t,T];\\
~~U_t =1.
\end{cases}
$$
\end{lemma}
 Remark 1 in \cite{Sam3} states that conditions \eqref{sam10} and \eqref{sam11} in Lemma \ref{Duality} can be simplified to
\begin{equation}\label{conds}
\eta_s \Psi_s^{\dagger}(e_j-X_{s-})>-1,
\end{equation}
for all $j\in\{1,\cdots,N\}$, without loss of generality. Note that if Assumption \ref{ass111} holds, we deduce condition (\ref{conds}) holds, furthermore, the result of Lemma \ref{Duality} holds.

\indent Lemma \ref{local} is Corollary 7.22 in Klebaner \cite{Klebaner}.
\begin{lemma}\label{local} Let $\{M(t); 0\leq t < \infty\}$ be a local martingale such that for all
$t$, $E[\sup\limits_{s\leq t}
|M(s)|] < \infty$. Then it is a martingale.
\end{lemma}

\section{An estimate of the solutions to anticipated BSDEs with Markov chain model}
\indent In order to make this
paper self-contained, we shall provide a proof of the existence and uniqueness of solutions of anticipated BSDEs with Markov chain noise by using the
fixed point theorem, rather than using Picard iterations as in Lu and Ren \cite{luren}.\\
\indent Consider the following anticipated BSDE on the Markov chain:
\begin{equation}\label{bfewwww10}
\left \{ \begin{array}{ll}
-dY_t = f(t,Y_t,Z_t,Y_{t+\delta(t)}, Z_{t+\zeta(t)} )dt- Z_t' dM_t,&t\in [ 0,T]; \\
\hskip.57cmY_{t}= \xi_t,&t\in[T,T+K];\\
\hskip.54cmZ_t=\eta_t,&t\in [ T,T+K].\end{array} \right.
\end{equation}Here $M$ is defined in (\ref{semimartingale}), $\delta(\cdot)$ and $\zeta(\cdot)$ are two $\mathbb{R}^+$-valued continuous functions defined on $[0,T]$ such that\\
${\bf (i)}$ there exists a constant $K\geq0$ such that
for any $s\in[0,T],$
$$ s+\delta(s)\leq T+K,\hskip2cm s+\zeta(s)\leq T+K;$$
${\bf (ii)}$ there exists a constant $L\geq0$ such that
for any $t\in[0,T]$ and a nonnegative and integrable function $g(\cdot)$,
$$ \int ^T _tg(s+\delta(s))ds\leq L\int ^{T+K} _tg(s)ds;$$
$$\int ^T _tg(s+\zeta(s))ds\leq L\int ^{T+K} _tg(s)ds.$$
\indent Assume that for any $s\in
[0,T],~f(s,\omega,y,z,\xi,\eta):\Omega\times \mathbb{R}\times
\mathbb{R}^{N}\times
L^{2}(\mathcal{F}_r;\mathbb{R})\times
H^{2}(\mathcal{F}_{r'};\mathbb{R}^{N})\longrightarrow
L^{2}(\mathcal{F}_{s},\mathbb{R}),$ where $r,r'\in[s,T+K]$, and $f$ satisfies the following conditions  \\
${\bf (H1)}$ There exist two constants $c_1,c_2>0,$ such that for any $s\in
[0,T],~y,y'\in\mathbb{R},$ $z,z'\in\mathbb{R}^{N},~\xi_{\cdot},\xi'_{\cdot}\in L^2_{\mathcal{F}}
(s,T+K;\mathbb{R}),~\eta_{\cdot},\eta'_{\cdot}\in
H^2 (s,T+K;\mathbb{R}^N),~r,\bar{r}\in[s,T+K],$ we have
$$\begin{array}{ll}
\hskip.46cm|f(s,y,z,\xi_r,\eta_{\bar{r}})-f(s,y',z',\xi'_r,\eta'_{\bar{r}})| \\[0.2cm]
\leq
c_1(|y-y'|+E^{\mathcal{F}_s}[|\xi_r-\xi'_r|])+c_2(\|z-z'\|_{X_s}+E^{\mathcal{F}_s}[\|\eta_{\bar{r}}-\eta'_{\bar{r}}\|_{X_s}]).
\end{array}
$$
${\bf (H2)}$ For each $(y,z,\xi,\eta)\in \mathbb{R} \times \mathbb{R}^N\times
L^{2}(\mathcal{F}_r;\mathbb{R})\times
H^{2}(\mathcal{F}_{r'};\mathbb{R}^{N})$, the process $(f(t,y,z,\xi,\eta))_{t\in [0,T]}$ is predictable, and $E[\int_{0} ^{T} |f(s,0,0,0,0)|^2ds] < \infty$.\\[0.2cm]
\indent Lu and Ren \cite{luren} proved the result of Theorem \ref{bft3000000} below. Here, we give an alternative proof.
\begin{theorem}\label{bft3000000} Suppose that $f$ satisfies $(H1)$ and $(H2),$ $\delta,\zeta$ satisfy $(i)$ and $(ii).$ Then for arbitrary given terminal
conditions $\xi_{\cdot}\in L^2_{\mathcal{F}}
(T,T+K;\mathbb{R}),$ $\eta_{\cdot}\in H^2
(T,T+K;\mathbb{R}^{N})$, the anticipated BSDE (\ref{bfewwww10}) has a
unique solution, i.e., there exists a unique pair of
stochastic processes $(Y.,Z.)\in
L^2_\mathcal{F}(0,T+K;\mathbb{R})\times H^2(0,T+K;\mathbb{R}^{N})$ satisfying equation (\ref{bfewwww10}). Moreover, this solution is
unique up to indistinguishability for $Y$ and equality $d\langle
X,X\rangle_t$ $\times\mathbb{P}$-a.s. for $Z$.\end{theorem}
 \noindent{\bf Proof.} Set $c:=\max\{c_1,c_2\}$. We fix
$\beta=16c^2(L+1),$ where $L$ is given in $(ii)$. Now we introduce a norm in the Banach space
$L^2_\mathcal{F}(0,T+K;\mathbb{R}):$\\[2mm]
\indent\indent\indent\indent\indent\indent$\|\nu_.\|_{L^2}=(E[\int ^{T+K} _0|\nu_s|^{2}e^{\beta s}ds])^{\frac{1}{2}}.$\\[2mm]
Define an equivalence class of $\varphi_.$ by $[\varphi_.]:=\{\psi_.;~E[\int_0^{T+K} \|\psi_t-\varphi_t\|^2_{X_s} ds]=0 \}$ and denote the factor space of equivalence classes of processes in $H^2(0,T+K;\mathbb{R}^{N})$ by $\hat{H}^2(0,T+K;\mathbb{R}^{N}):=\{[\varphi_.];~\varphi_.\in H^2(0,T+K;\mathbb{R}^{N}) \}.$ Then $\hat{H}^2(0,T+K;\mathbb{R}^{N})$ is a Banach space with the norm \\[2mm]
\indent\indent\indent\indent\indent\indent$\|\mu_.\|_{\hat{H}^2}=(E [ \int_0^{T+K}\|\mu_s\|^2_{X_s}e^{\beta s} ds ])^{\frac{1}{2}}.$\\[2mm]
 Set
$$ \left \{ \begin{array}{ll}Y_t =\xi_T +  \int ^T _t f(s,y_s,z_s,y_{s+\delta(s)},z_{s+\zeta(s)} )ds   - \int^{T} _{t} Z_s'dM_s,&t\in[0,T]; \\
Y_{t} = \xi_t,&t\in[T,T+K];\\
Z_t=\eta_t,&t\in[T,T+K].\end{array} \right.$$ By Lemma \ref{existence}, we know for any
$(y_{\cdot},z_{\cdot}) \in L^2 _\mathcal{F}(0,T+K;\mathbb{R})\times H^2(0,T+K;\mathbb{R}^{N})$,
the above equation has a solution $(Y_{\cdot},Z_{\cdot})\in L^2 _\mathcal{F}(0,T+K;\mathbb{R})\times H^2(0,T+K;$ $\mathbb{R}^{N})$, moreover, this solution is
unique up to indistinguishability for $Y_{\cdot}$ and equality $d\langle
X,X\rangle_t\times\mathbb{P}$-a.s. for $Z_{\cdot}$. That is, this solution is
unique up to indistinguishability for $(Y_{\cdot},Z_{\cdot})\in L^2 _\mathcal{F}(0,T+K;\mathbb{R})\times \hat{H}^2(0,T+K;\mathbb{R}^{N})$. Define a mapping $h
:L^2 _\mathcal{F}(0,T+K;\mathbb{R})\times \hat{H}^2(0,T+K;\mathbb{R}^{N})\longrightarrow L^2 _\mathcal{F}(0,T+K;$ $\mathbb{R})\times \hat{H}^2(0,T+K;\mathbb{R}^{N})$ such
that $h[(y_{\cdot},z_{\cdot})]=(Y_{\cdot},Z_{\cdot})$. Now we prove
that $h$ is a contraction mapping under the norm
$\|\cdot\|_{L^2}+\|\cdot\|_{\hat{H}^2}$. For two arbitrary elements
$(y_{\cdot},z_{\cdot})$ and $(y'_{\cdot},z'_{\cdot})$ in
$L^2 _\mathcal{F}(0,T+K;\mathbb{R})\times \hat{H}^2(0,T+K;\mathbb{R}^{N})$
set $(Y_{\cdot},Z_{\cdot})=h[(y_{\cdot},z_{\cdot})]$ and
$(Y'_{\cdot},Z'_{\cdot})=h[(y'_{\cdot},z'_{\cdot})]$. Denote their
differences by
$(\hat{y}_{\cdot},\hat{z}_{\cdot})=((y-y')_{\cdot},(z-z')_{\cdot})$ and $(\hat{Y}_{\cdot},\hat{Z}_{\cdot})=((Y-Y')_{\cdot},(Z-Z')_{\cdot}).$
Applying Product Rule for Semimartingales  in \cite{elliott} to $|\hat{Y}_t|$, we have
\begin{align*}
|\hat{Y}_t|^2
& = -2 \int_t^T \hat{Y}_{s-} d\hat{Y}_{s} - \sum_{t \leq s \leq T} \Delta \hat{Y}_s\Delta  \hat{Y}_s \\
&= -2 \int_t^T \hat{Y}_s (f(s,y_s,z_s,y_{s+\delta(s)},z_{s+\zeta(s)})-f(s,y'_s,z'_s,y'_{s+\delta(s)},z'_{s+\zeta(s)})) ds \\
& \quad -2 \int_t^T  \hat{Y}_{s-} (\hat{Z}_s)' dM_s - \sum_{t \leq s \leq T} \Delta \hat{Y}_{s}\Delta \hat{Y}_{s}.
\end{align*}
Also
\begin{align*}
&\sum_{t\leq s \leq T}  \Delta  \hat{Y}_{s} \Delta  \hat{Y}_{s}
= \sum_{t\leq s \leq T}( (\hat{Z}_s)' \Delta X_s)(  ( \hat{Z}_s)'\Delta X_s ) = \sum_{t\leq s \leq T} (\hat{Z}_s )' \Delta X_s \Delta X_s' \hat{Z}_s \\
&=  \int_t^T ( \hat{Z}_s )'  (dL_s + d\lel X,X\rir_s) \hat{Z}_s=  \int_t^T (\hat{Z}_s)' dL_s\hat{Z}_s + \int_t^T \|\hat{Z}_s \|_{X_s}^2 ds.
\end{align*}
Applying It\^{o}'s formula to $e^{\beta s} |\hat{Y}_s|^2$ for $s\in[0,T]$ and then taking the expectation:
\begin{align*}
&E [  |\hat{Y}_0|^2 ]+ E [ \int_0^T \beta |\hat{Y}_s|^2 e^{\beta s} ds ]  + E [ \int_0^T\|\hat{Z}_s\|^2_{X_s}e^{\beta s} ds ] \\
& = 2 E [ \int_0^T \hat{Y}_s (f(s,y_s,z_s,y_{s+\delta(s)},z_{s+\zeta(s)})-f(s,y'_s,z'_s,y'_{s+\delta(s)},z'_{s+\zeta(s)})) e^{\beta s}ds ]\\
&\leq  E [ \int_0^T (\frac{\beta}{2}|\hat{Y}_s|^2+\frac{2}{\beta} |f(s,y_s,z_s,y_{s+\delta(s)},z_{s+\zeta(s)})-f(s,y'_s,z'_s,y'_{s+\delta(s)},z'_{s+\zeta(s)})|^2)e^{\beta s} ds ].
\end{align*}
Since $\delta(s),\zeta(s)$ satisfy $(ii)$ and $f$ satisfies $(H1)$,
by the Fubini Theorem we have
$$\begin{array}{lll}
E[\int^T_0(\dfrac{\beta}{2}|\hat{Y}_s|^{2}+\|\hat{Z}_s \|_{X_s}^{2})e^{\beta
s}ds]\\[0.3cm]
\leq\dfrac{2c^2}{\beta}E[\int^{T}_0(|\hat{y}_s|+\|\hat{z}_s\|_{X_s}+E^{\mathcal{F}_s}[~|\hat{y}_{s+\delta(s)}|+\|\hat{z}_{s+\zeta(s)}\|_{X_s}])^2e^{\beta
s}ds]\\[0.3cm]
\leq\dfrac{8c^2}{\beta}E[\int^{T}_0(|\hat{y}_s|^2+\|\hat{z}_s\|_{X_s}^2+|\hat{y}_{s+\delta(s)}|^2+\|\hat{z}_{s+\zeta(s)}\|_{X_s}^2)e^{\beta
s}ds]\\[0.3cm]
\leq\dfrac{8c^2(L+1)}{\beta}E[\int^{T+K}_0(|\hat{y}_s|^{2}+\|\hat{z}_s\|_{X_s}^{2})e^{\beta
s}ds].
\end{array}$$
Note $\beta=16c^2(L+1),$ therefore
$$E[\int ^{T+K} _0(|\hat{Y}_s|^{2}+\|\hat{Z}_s\|_{X_s}^{2})e^{\beta
s}ds]\leq\frac{1}{2}E[\int
^{T+K}_0(|\hat{y}_s|^{2}+\|\hat{z}_s\|_{X_s})e^{\beta s}ds],
$$
or $$\|\hat{Y}_{\cdot}\|_{L^2}+\|\hat{Z}_{\cdot}\|_{\hat{H}^2}
\leq\frac{1}{\sqrt{2}}(\|\hat{y}_{\cdot}\|_{L^2}+\|\hat{z}_{\cdot}\|_{\hat{H}^2}).$$
Consequently $h$ is a strict contraction mapping on
$L^2_\mathcal{F}(0,T+K;\mathbb{R})\times \hat{H}^2(0,T+K;\mathbb{R}^{N}).$
It follows by the Fixed Point Theorem that the anticipated BSDE (\ref{bfewwww10})
has a unique solution $(Y_{\cdot},Z_{\cdot})\in
L^2_\mathcal{F}(0,T+K;\mathbb{R})\times \hat{H}^2(0,T+K;\mathbb{R}^{N}).$ That is, the solution $(Y_{\cdot},Z_{\cdot})\in
L^2_\mathcal{F}(0,T+K;\mathbb{R})\times H^2(0,T+K;\mathbb{R}^{N})$ is
unique up to indistinguishability for $Y$ and equality $d\langle
X,X\rangle_t\times\mathbb{P}$-a.s. for $Z$.
$\mbox{}\hfill\Box$\\[2mm]
\indent Our method allows us to find an estimate of the solution to equation
(\ref{bfewwww10}).
\begin{prop} \label{prop}Assume that $f$ satisfies $(H1)$ and $(H2),$ $\delta$ and
$\zeta$ satisfy $(i)$ and $(ii)$. Then there exists a
constant $C>0$ depending only on $c_1,c_2$ in $(H1),$ $L$ in $(ii)$, $K$ and
$T $, such that for each $\xi_{\cdot}\in L^2_{\mathcal{F}}
(T,T+K;\mathbb{R}),~\eta_{\cdot}\in H^2
(T,T+K;\mathbb{R}^{N}),$ the solution $(Y.,Z.)$ to the
anticipated BSDE (\ref{bfewwww10}) satisfies
\begin{equation}\label{bfe1000}\begin{array}{lll}
E[\sup \limits_{0\leq s\leq T }|Y_s|^{2}+\int^T_0\|Z_s\|_{X_s}^{2}ds]\\[0.3cm]
 \leq
CE[|\xi_{T}|^{2}+\int^{T+K}_T(|\xi_s|^{2}+\|\eta_s\|_{X_s}^{2})ds+(\int^T_0|f(s,0,0,0,0)|ds)^{2}].
\end{array}\end{equation}
\end{prop}
\noindent{\bf Proof.} Set $c=:\max\{c_1,c_2\}$. let $\beta >0$ be an arbitrary constant. Using It$\hat{\text{o}}$'s formula for $e^{\beta t}|Y_t|^2$, we deduce
$$\begin{array}{ll}
E [  |Y_0|^2 ]+ E [ \int_0^T \beta |Y_s|^2 e^{\beta s} ds ]  + E[ \int_0^T e^{\beta s}\|Z_s\|^2_{X_s} ds ] \\[3mm]
= E [ e^{\beta T} |\xi_T|^2 ] + 2 E [ \int_0^T e^{\beta s} Y_s f(s, Y_s, Z_s,Y_{s+\delta(s)}, Z_{s+\zeta(s)} ) ds ]\\[3mm]
\leq  E [e^{\beta T}|\xi_T|^2] + 2 E[ \int_0^T e^{\beta s} |Y_s|\cdot|f(s,0,0,0,0)|ds]\\[3mm]
+2 E[ \int_0^T e^{\beta s}|Y_s|\cdot |f(s, Y_s, Z_s,Y_{s+\delta(s)}, Z_{s+\zeta(s)} )-f(s,0,0,0,0)|ds] \\[3mm]
\leq  E [e^{\beta T}|\xi_T|^2] + 2 E[ \int_0^T e^{\beta s}|Y_s| \cdot|f(s, 0, 0,0,0)| ds]\\[3mm]
+2c E[ \int_0^T e^{\beta s}|Y_s|( |Y_s|+E^{\mathcal{F}_s}[|Y_{s+\delta(s)}|] + \|Z_s\|_{X_s}+E^{\mathcal{F}_s}[\|Z_{s+\zeta(s)} \|_{X_s}])ds] \\[3mm]
 \leq E [e^{\beta T}|\xi_T|^2]  +2 E [ \sup\limits_{s\in[0,T]}e^{\frac{1}{2}\beta s}|Y_s|\cdot \int_0^Te^{\frac{1}{2}\beta s} |f(s,0,0,0,0)| ds]\\[3mm]
  + (3c+3c^2+3c^2L) E[\int_0^T e^{\beta s}|Y_s|^2  ds]  + cE[\int_0^T e^{\beta s}|Y_{s+\delta(s)}|^2  ds]\\[3mm]
 + \dfrac{1}{3} E [ \int_0^T e^{\beta s}\|Z_s\|^2_{X_s} ds]+ \dfrac{1}{3L} E [ \int_0^T e^{\beta s}\|Z_{s+\zeta(s)}\|^2_{X_s} ds]\\[3mm]
\leq E [e^{\beta T}|\xi_T|^2] +\alpha E [ \sup\limits_{s\in[0,T]}e^{\beta s}|Y_s|^2]+\dfrac{1}{\alpha}E[ (\int_0^Te^{\frac{1}{2}\beta s} |f(s, 0,0,0, 0)| ds)^2] \\[3mm]
 + (3c+3c^2+3c^2L+cL) E[\int_0^{T+K} e^{\beta s}|Y_s|^2  ds] + \dfrac{2}{3} E [ \int_0^{T+K} e^{\beta s}\|Z_s\|^2_{X_s} ds],
\end{array}
$$
where $\alpha>0$ is also an arbitrary constant. Set $\beta=3c+3c^2+3c^2L+cL+1$, we obtain
\begin{align}\label{apriory21}
&
E[ \int_0^T |Y_s|^2 e^{\beta s} ds] + \frac{1}{3} E [ \int_0^T e^{\beta s}\|Z_s\|^2_{X_s} ds]  \\
\nonumber&\leq E [e^{\beta T}|\xi_T|^2] +\alpha E [ \sup\limits_{s\in[0,T]}e^{\beta s}|Y_s|^2]+\dfrac{1}{\alpha}E[ (\int_0^Te^{\frac{1}{2}\beta s} |f(s, 0,0,0, 0)| ds)^2]\\
\nonumber& + (3c+3c^2+3c^2L+cL) E[\int_T^{T+K} e^{\beta s}|\xi_s|^2  ds] + \dfrac{2}{3} E [ \int_T^{T+K} e^{\beta s}\|\eta_s\|^2_{X_s} ds].
\end{align}
Using Doob's inequality and Lemma \ref{Z2}, we know
\begin{align}\label{marsup}
 \nonumber &E [ \sup_{0\leq t \leq T}|\int_t^T Z_s' dM_s|^2 ] =E [ \sup_{0\leq t \leq T}|\int_0^T Z_s' dM_s-\int_0^t Z_s' dM_s|^2]\\
 \nonumber&\leq 2 E [|\int_0^T Z_s' dM_s|^2+ \sup_{0\leq t \leq T}|\int_0^t Z_s' dM_s|^2] \leq 10 E [ |\int_0^T Z'_sdM_s|^2]\\
&= 10 E[ \int_0^T \|Z_s\|^2_{X_s}ds]\leq10 E[ \int_0^T e^{\beta s}\|Z_s\|^2_{X_s}ds].
\end{align}
Because $ Y_t = \xi + \int_t^T f(s, Y_s, Z_s,Y_{s+\delta(s)}, Z_{s+\zeta(s)}) ds-\int_t^T  Z_s' dM_s,~  0 \leq t \leq T$, by (\ref{marsup}) we have
\begin{align}\label{supp}
\nonumber&E[\sup\limits_{0\leq t \leq T}|Y_t|^2]\\
\nonumber&\leq E[3|\xi_T|^2 + 3(\int_0^T| f(s, Y_s, Z_s,Y_{s+\delta(s)}, Z_{s+\zeta(s)}) |ds)^2 + 3 \sup\limits_{0\leq t \leq T}|\int_t^T Z_s' dM_s|^2] \\
\nonumber&\leq 3E[|\xi_T|^2]+ 30E[ \int_0^T e^{\beta s}\|Z_s\|^2_{X_s}ds] \\
\nonumber&+ 3E[( \int_0^T (|f(s,0,0,0,0)|+c|Y_s|+c\|Z_s\|_{X_s}+c|Y_{s+\delta(s)}|+c\|Z_{s+\zeta(s)}\|_{X_s})ds )^2 ] \\
\nonumber& \leq 3E[|\xi_T|^2]+ 30E[ \int_0^T e^{\beta s}\|Z_s\|^2_{X_s}ds]+ 15E[( \int_0^T |f(s,0,0,0,0)|ds)^2] \\
 \nonumber&+ 15Tc^2E[\int_0^T( |Y_s|^2+\|Z_s\|^2_{X_s}+|Y_{s+\delta(s)}|^2+\|Z_{s+\zeta(s)}\|_{X_s}^2)ds] \\
\nonumber&\leq 3E[e^{\beta T}|\xi_T|^2]+15E[( \int_0^T |f(s,0,0,0,0)|ds)^2] \\
\nonumber&+ 15(2 +Tc^2+Tc^2L)E[\int_0^{T}e^{\beta s} (|Y_s|^2+\|Z_s\|^2_{X_s})ds]\\
&+ 15Tc^2LE[\int_T^{T+K}e^{\beta s} (|\xi_s|^2+\|\eta_s\|^2_{X_s})ds].
\end{align}
Set $\alpha=\dfrac{1}{90(2+Tc^2+Tc^2L)}$. Then by (\ref{apriory21}) and (\ref{supp}), we deduce there exists a constant $C>0$ depending on $T,c,L$ and $K$ such that (\ref{bfe1000}) holds.$\mbox{}\hfill\Box$
\section{Duality between SDDEs and Anticipated BSDEs on Markov chains}
  It is well known that there is perfect duality
 between SDEs and BSDEs (see El Karoui, Peng, and Quenez \cite{EPQ1997}). Cohen, Elliott \cite{Sam2} and \cite{Sam3} showed duality between SDEs and BSDEs driven by Markov chains. In \cite{PengYang} Peng and Yang considered duality between SDDEs and anticipated
 BSDEs. We now establish duality between SDDEs and anticipated
 BSDEs with Markov chain noise.
 \begin{lemma}\label{normcon}
For any $B \in \mathbb{R}^{N\times N}$,
$$ ~~~~\|B\|_{X_t}^2 \leq 3m \|B\|_{N\times N}^2, ~~\text{ for any }t\in[0,T],$$
where $m>0$ is the bound of $\|A_t\|_{N\times N}$, for any $t\in[0,T]$.
\end{lemma}
 \noindent{\bf Proof.} Write $B=(B_1,B_2,\ldots,B_N)$, where $B_i\in \mathbb{R}^{N}$, for any $1\leq i\leq N$. Then $\|B\|_{N\times N}^2=\sum^{N}_{i=1}|B_i|_N^2.$
 Noticing that for any $1\leq i\leq N$, $B_i'\Psi_tB_i\in\mathbb{R}$, we obtain for any $t\in[0,T]$,
 \begin{align*}
\|B\|_{X_t}^2 &=Tr((B_1,B_2,\ldots,B_N)'\Psi_t(B_1,B_2,\ldots,B_N))\\[2mm]
& = Tr((B_1'\Psi_tB_1,B_2'\Psi_tB_2,\ldots,B_N'\Psi_tB_N))\\[2mm]
&=\sum^{N}_{i=1}B_i'\Psi_tB_i=\sum^{N}_{i=1}\|B_i\|_{X_t}^2.
\end{align*}
By Lemma \ref{normbound} we have
$\|B\|_{X_t}^2\leq 3m\sum^{N}_{i=1}|B_i|_{N}^2=3m\|B\|_{N\times N}^2,$ for any $t\in[0,T].$$\mbox{}\hfill\Box$
 \begin{ass}\label{ass0}
Assume there exists a constant $l>0$ such that for any $t\in[0,T]$, $\|\Psi_t^{\dagger}\|_{N\times N}^2\leq l,$ where $\Psi$ is given in \eqref{Psi}.
\end{ass}
 \begin{lemma}\label{Zequality} Suppose that Assumption \ref{ass0} holds, $f$ satisfies $(H1),(H2)$ and $\delta,\zeta$ satisfy $(i)$ and $(ii).$ Then for any $\xi_{\cdot}\in L^2_{\mathcal{F}}
(T,T+K;\mathbb{R}),\eta_{\cdot}\in H^2
(T,T+K;$ $\mathbb{R}^{N})$, the solution $Z.\in H^2(0,T+K;\mathbb{R}^{N})$ of the anticipated BSDE (\ref{bfewwww10}) satisfies $Z.=(\Psi\Psi^{\dag}Z).$, $d\langle
X,X\rangle_t$ $\times\mathbb{P}$-a.s.\end{lemma}
 \noindent{\bf Proof.} Set $c:=\max\{c_1,c_2\}$. By the proof of Theorem \ref{bft3000000},We know there exists a sequence of
$\{(y^{(n)}_{\cdot},z^{(n)}_{\cdot});~n\in\mathbb{N}\}\subseteq  L^2 _\mathcal{F}(0,T+K;\mathbb{R})\times H^2(0,T+K;$ $\mathbb{R}^{N})$ satisfying for any $n\in\mathbb{N},$
$$ \left \{ \begin{array}{ll}y^{(n+1)}_t =\xi_T +  \int ^T _t f(s,y^{(n+1)}_s,z^{(n+1)}_s,y^{(n)}_{s+\delta(s)},z^{(n)}_{s+\zeta(s)} )ds   - \int^{T} _{t} (z^{(n+1)}_s)'dM_s,~t\in[0,T]; \\
y^{(n+1)}_{t} = \xi_t,~~~~t\in[T,T+K];\\
z^{(n+1)}_t=\eta_t,~~~~t\in[T,T+K].\end{array} \right.$$
Then $$E[\int ^{T+K} _0(|y^{(n)}_s-Y_s|^{2}+\|z_s^{(n)}-Z_s\|_{X_s}^{2})e^{\beta
s}ds]\rightarrow0, ~~~~~\mbox{as }~~n\rightarrow\infty,
$$
where $(Y_{\cdot},Z_{\cdot})\in L^2 _\mathcal{F}(0,T+K;\mathbb{R})\times H^2(0,T+K;\mathbb{R}^{N})$ is the solution of the anticipated BSDE (\ref{bfewwww10}). Thus, $E[\int ^{T+K} _0\|z_s^{(n)}-Z_s\|_{X_s}^{2}
ds]\rightarrow0$ as $n\rightarrow\infty$. By Lemma \ref{sam35}, we have for any $n\in\mathbb{N}$, $E[\int ^{T+K} _0\|z_s^{(n)}-\Psi_s\Psi_s^\dag z_s^{(n)}\|_{X_s}^{2}ds]=0.$ Noting $\Psi_t=\diag(A_tX_t)- \diag(X_t)A'_t - A_t \diag(X_t)$ given in (\ref{Psi}), by Lemma \ref{normcon} we obtain for any $t\in[0,T]$,
\begin{align*}
\|\Psi_t\|^2_{X_t}&\leq 3m \|\diag(A_tX_t)- \diag(X_t)A'_t - A_t \diag(X_t)\|_{N\times N}^2\\
& \leq  3m (|A_tX_t|_{ N}+|X_t|_N\cdot\|A_t\|_{N\times N} +\| A_t\|_{N\times N}\cdot|X_t|_N )^2\\
& \leq  3m(\|A_t\|_{N\times N}\cdot|X_t|_{ N}+|X_t|_N\cdot\|A_t\|_{N\times N} +\| A_t\|_{N\times N}\cdot|X_t|_N )^2\\
& \leq  27m\|A_t\|_{N\times N}^2 \leq 27m^3.
\end{align*}
Hence, by Assumption \ref{ass0} and Lemma \ref{normcon}, we deduce
\begin{align*}
&E[\int ^{T+K} _0\|\Psi_s\Psi_s^\dag  z_s^{(n)}-\Psi_s\Psi_s^\dag  Z_s\|_{X_s}^{2}ds]\\
&\leq E[\int ^{T+K} _0\|\Psi_s\|_{X_s}^2\cdot\|\Psi^\dag_s\|_{X_s}^2 \cdot \|z_s^{(n)}- Z_s\|_{X_s}^{2}ds]\\
&\leq27m^3l E[\int ^{T+K} _0 \|z_s^{(n)}- Z_s\|_{X_s}^{2}ds]\rightarrow0,~~~~~~\mbox{as }~~n\rightarrow\infty.
\end{align*}
Therefore,
\begin{align*}
~~&E[\int ^{T+K} _0\|Z_s-\Psi_s\Psi_s^\dag  Z_s\|_{X_s}^{2}ds]=\lim_{n\rightarrow\infty}E[\int ^{T+K} _0\|Z_s-\Psi_s\Psi_s^\dag  Z_s\|_{X_s}^{2}ds]\\
&\leq 3\lim_{n\rightarrow\infty}E[\int ^{T+K} _0\|Z_s-z_s^{(n)}\|_{X_s}^{2}ds]+3\lim_{n\rightarrow\infty}E[\int ^{T+K} _0\|  z_s^{(n)}-\Psi_s\Psi_s^\dag  z^{(n)}_s\|_{X_s}^{2}ds]\\
&+3\lim_{n\rightarrow\infty}E[\int ^{T+K} _0\|\Psi_s\Psi_s^\dag  z_s^{(n)}-\Psi_s\Psi_s^\dag  Z_s\|_{X_s}^{2}ds]=0.~~~~~~~~~~~~~~~~~~~~~~~~~~~~~~~\mbox{}\hfill\Box
\end{align*}
 \begin{theorem}\label{bft1} Suppose $\theta>0$ is a given constant, $a_.,\mu_.\in
L^{2}_\mathcal{F}(t_0-\theta,T+\theta;\mathbb{R}),$ $\varphi_.\in
L^{2}_\mathcal{F}(t_0,T;\mathbb{R}),~b_.\in
 L^2_\mathcal{F}(t_0-\theta,T+\theta;\mathbb{R}^{1\times N}),$ and  moreover, there is a constant $\gamma>0$ such that $|a_s|\leq\gamma$, $|b_s|_N\leq\gamma$ and $|\mu_s|\leq\gamma$ for any $s\in[t_0-\theta,T+\theta]$. Then for all $ U_.\in
L^2_{\mathcal{F}} (T,T+\theta;\mathbb{R}),$ the solution $Y_.$ to
 anticipated BSDE with Markov chain noise
\begin{equation*}
\left \{ \begin{array}{ll}
-dY_{t}=(a_{t}Y_{t}+\mu_{t}E^{\mathcal{F}_{t}}[Y_{t+\theta
}]+b_{t}Z_{t}+\varphi_t)dt-Z'_{t}dM_{t},~t\in[t_0,T]; \\
\hskip.56cmY_{t}=U_t,~~~~t\in [T,T+\theta].
\end{array}
\right. \end{equation*}
can be given by the closed formula:
$$
Y_{t}
=E^{\mathcal{F}_{t}}[\hat{X}_TU_T+\int^T_{t}\hat{X}_s\varphi_sds+\int^{T+\theta}_{T}\mu_{s-\theta}\hat{X}_{s-\theta}U_sds],
$$
for any $ t \in [t_0,T]$, a.s., where $\hat{X}_s$ is the solution to SDDE with Markov chain
\begin{equation*}
\left \{ \begin{array}{ll}
d\hat{X}_{s}=(a_{s}\hat{X}_{s}+\mu_{s-\theta}\hat{X}_{s-\theta})ds+\hat{X}_{s-}b_{s-}
(\Psi_s^\dagger)'dM_{s},~s\in [ t,T+\theta];\\
\hskip.245cm\hat{X}_{t}
=1,\\
\hskip.245cm\hat{X}_s =0,~~~~ s\in [ t-\theta,t).
\end{array}
\right.\end{equation*}\end{theorem}
 \noindent {\bf Proof.} By Theorem 3.1 of Mao \cite{mao}, we have there exists a unique RCLL adapted solution $hat{X}$ of the above SDDE.
 By (\ref{semimartingale}), $[M,M]_t=[X,X]_t=\lel X,X \rir_t +L_t$ and $d \lel X,X\rir_t = \Psi_t dt$. By Definition \ref{defMoore} and Lemma \ref{Zequality}, $Z_t'=(\Psi_t\Psi_t^{\dagger} Z_t)'=Z_t'(\Psi_t\Psi^{\dagger}_t)'=Z_t'\Psi_t\Psi^{\dagger}_t$ for $t\in[t_0,T]$. Applying It\^{o}'s formula to
$\hat{X}_sY_s$ for $s\in[t,T],$ we derive
\begin{align*}
&d(\hat{X}_sY_s)\\
 & = \hat{X}_{s-} dY_s + Y_{s-}d\hat{X}_s + d[\hat{X},Y]_s \\[.5mm]
& = - \hat{X}_sa_{s}Y_{s}ds- \hat{X}_s\mu_{s}E^{\mathcal{F}_{s}}[Y_{s+\theta
}]ds- \hat{X}_sb_{s}Z_{s}ds- \hat{X}_s\varphi_sds+ \hat{X}_{s-}Z'_{s} dM_s  + Y_s \hat{X}_s a_s ds\\[.5mm]
&+ Y_s\mu_{s-\theta}\hat{X}_{s-\theta}ds+ Y_{s-}\hat{X}_{s-}b_{s-}(\Psi_s^\dagger)'dM_{s}
+ Z'_s \Delta M_s \hat{X}_{s-}b_{s-} (\Psi_s^{\dagger})' \Delta M_s\\[.5mm]
& =- \hat{X}_s\mu_{s}E^{\mathcal{F}_{s}}[Y_{s+\theta
}]ds- \hat{X}_sb_{s}Z_{s}ds- \hat{X}_s\varphi_sds+ \hat{X}_{s-}Z'_{s} dM_s  + Y_s\mu_{s-\theta}\hat{X}_{s-\theta}ds \\[.5mm]
&+ Y_{s-}\hat{X}_{s-}b_{s-}(\Psi_s^\dagger)'dM_{s}
+  Z'_s \Delta M_s \Delta M'_s \Psi_s^{\dagger} \hat{X}_{s-}b'_{s-} \\[.5mm]
& =- \hat{X}_s\mu_{s}E^{\mathcal{F}_{s}}[Y_{s+\theta
}]ds- \hat{X}_sb_{s}Z_{s}ds- \hat{X}_s\varphi_sds+ \hat{X}_{s-}Z'_{s} dM_s  + Y_s\mu_{s-\theta}\hat{X}_{s-\theta}ds \\[.5mm]
&+ Y_{s-}\hat{X}_{s-}b_{s-}(\Psi_s^\dagger)'dM_{s}
+  Z'_s d[M ,M]_s\Psi_s^{\dagger} \hat{X}_{s-}b'_{s-}\\[.5mm]
& = - \hat{X}_s\mu_{s}E^{\mathcal{F}_{s}}[Y_{s+\theta
}]ds- \hat{X}_sb_{s}Z_{s}ds- \hat{X}_s\varphi_sds+ \hat{X}_{s-}Z'_{s} dM_s + Y_s\mu_{s-\theta}\hat{X}_{s-\theta}ds \\[.5mm]
& + Y_{s-}\hat{X}_{s-}b_{s-}(\Psi_s^\dagger)'dM_{s} + Z'_s\Psi_s \Psi_s^{\dagger}\hat{X}_sb'_s ds  + Z'_sdL_s \Psi_s^{\dagger} \hat{X}_{s-}b'_{s-}\\[.5mm]
&= - \hat{X}_s\mu_{s}E^{\mathcal{F}_{s}}[Y_{s+\theta
}]ds- \hat{X}_sb_{s}Z_{s}ds -\hat{X}_s\varphi_sds+ \hat{X}_{s-}Z'_{s} dM_s  + Y_s\mu_{s-\theta}\hat{X}_{s-\theta}ds\\[.5mm]
&+ Y_{s-}\hat{X}_{s-}b_{s-}(\Psi_s^\dagger)'dM_{s}+ Z'_s\hat{X}_sb'_s ds + Z'_sdL_s \Psi_s^{\dagger} \hat{X}_{s-}b'_{s-}\\[.5mm]
&= - \hat{X}_s\mu_{s}E^{\mathcal{F}_{s}}[Y_{s+\theta
}]ds- \hat{X}_s\varphi_sds+ \hat{X}_{s-}Z'_{s} dM_s  + Y_s\mu_{s-\theta}\hat{X}_{s-\theta}ds\\[.5mm]
&+ Y_{s-}\hat{X}_{s-}b_{s-}(\Psi_s^\dagger)'dM_{s}+ Z'_sdL_s \Psi_s^{\dagger} \hat{X}_{s-}b'_{s-}.
\end{align*}
Then for any $s\in[t,T]$, we obtain
\begin{align*}
\hat{X}_sY_s-Y_t+\int^s_t \hat{X}_r\mu_{r}E^{\mathcal{F}_{r}}[Y_{r+\theta
}]dr+\int^s_t\hat{X}_r\varphi_rdr-\int^s_tY_r\mu_{r-\theta}\hat{X}_{r-\theta}dr
= \tilde{L}_s
\end{align*}
for some local martingale $\tilde{L}$. Thus by  H\"{o}lder's inequality, noting $\hat{X}_s =0$ for any $ s\in [ t-\theta,t)$, we know for any $T'\in[t,T]$,
\begin{align*}
&E[\sup_{s\in[t,T']}|\tilde{L}_s|]\\
&\leq E[\sup_{s\in[t,T']}|\hat{X}_sY_s|]+E[|Y_t|]+\gamma E[\sup_{s\in[t,T']}\int^s_t| \hat{X}_rE^{\mathcal{F}_{r}}[Y_{r+\theta
}]|dr]\\
&+ E[\sup_{s\in[t,T']}\int^s_t|\hat{X}_r\varphi_r|dr] +\gamma E[\sup_{s\in[t,T']}\int^s_t|Y_r\hat{X}_{r-\theta}|dr]
\end{align*}
\begin{align*}
&\leq \frac{1}{2}E[\sup_{s\in[t,T]}|\hat{X}_s|^2+\sup_{s\in[t,T]}|Y_s|^2]+E[|Y_t|]+\gamma E[\int^{T}_t| \hat{X}_rY_{r+\theta
}|dr]\\
&+ E[\int^{T}_t|\hat{X}_r\varphi_r|dr] +\gamma E[\int^{T}_t|Y_r\hat{X}_{r-\theta}|dr]\\
&\leq \frac{1}{2}E[\sup_{s\in[t,T]}|\hat{X}_s|^2+\sup_{s\in[t,T]}|Y_s|^2]+E[|Y_t|]+\gamma( E[\int^{T}_t|\hat{X}_r|^2dr])^{\frac{1}{2}}( E[\int^{T}_t|\varphi_r|^2dr])^{\frac{1}{2}}\\
&+\gamma (E[\int^{T}_t| \hat{X}_r|^2dr])^{\frac{1}{2}}(E[\int^{T}_t| Y_{r+\theta
}|^2dr])^{\frac{1}{2}}+\gamma (E[\int^{T}_t|Y_r|^2dr])^{\frac{1}{2}}(E[\int^{T}_t|\hat{X}_{r-\theta}|^2dr])^{\frac{1}{2}}\\
&\leq \frac{1}{2}E[\sup_{s\in[t,T]}|\hat{X}_s|^2+\sup_{s\in[t,T]}|Y_s|^2]+E[|Y_t|]+\gamma( E[\int^{T}_t|\hat{X}_r|^2dr])^{\frac{1}{2}}( E[\int^{T}_t|\varphi_r|^2dr])^{\frac{1}{2}}\\
&+\gamma (E[\int^{T}_t| \hat{X}_r|^2dr])^{\frac{1}{2}}(E[\int^{T+\theta
}_{t+\theta
}| Y_{r
}|^2dr])^{\frac{1}{2}}+\gamma (E[\int^{T}_t|Y_r|^2dr])^{\frac{1}{2}}(E[\int^{T-\theta}_t|\hat{X}_{r}|^2dr])^{\frac{1}{2}}.
\end{align*}
\begin{lemma}\label{Xbounded}
$$ E[\int^{T}_t|\hat{X}_{r}|^2dr]<+\infty,\text{ moreover, }E[\sup_{s\in[t,T]}|\hat{X}_s|^2]<+\infty.$$
\end{lemma}
\begin{proof} By Lemma \ref{Psibound}, we know there exists a constant $\rho>0$ such that $|\Psi_s^\dagger|_{N\times N}\leq \rho$ for $s\in[t_0-\theta,T+\theta]$. Since $\hat{X}_{s}=0$ for $s\in[t-\theta,t)$, by Lemma \ref{normbound} we have for $s\in[t,T+\theta]$,
\begin{align*}
&E[|\hat{X}_{s}|^2]\\
& \leq4(1+E[|\int^s_ta_{r}\hat{X}_{r}dr|^2]+E[|\int^s_t\mu_{r-\theta}\hat{X}_{r-\theta}dr|^2]+E[|\int^s_t\hat{X}_{r-}b_{r-}
(\Psi_r^\dagger)'dM_{r}|^2])\\
& \leq4+4\gamma^2(s-t)E[\int^s_t|\hat{X}_{r}|^2 dr]+4\gamma^2(s-t)E[\int^{s-\theta}_{t}|\hat{X}_{r}|^2 dr]\\
&+4E[\int^s_t\|\hat{X}_{r}b_{r}
(\Psi_r^\dagger)'\|^2_{X_r }dr]\\
& \leq4+8\gamma^2(T+\theta-t)E[\int^s_t|\hat{X}_{r}|^2 dr]
+12m^2E[\int^s_t|\hat{X}_{r}b_{r}
(\Psi_r^\dagger)'|^2_{N }dr]\\
& \leq4+8\gamma^2(T+\theta-t)E[\int^s_t|\hat{X}_{r}|^2 dr]
+12m^2E[\int^s_t|\hat{X}_{r}|^2\cdot|b_{r}|_N^2\cdot|
\Psi_r^\dagger|^2_{N\times N}dr]\\
& \leq4+8\gamma^2(T+\theta-t)E[\int^s_t|\hat{X}_{r}|^2 dr]
+12m^2\gamma^2\rho^2E[\int^s_t|\hat{X}_{r}|^2dr].
\end{align*}
By Gr\"{o}nwall's inequality, we derive for $s\in[t,T+\theta]$,
$$E[|\hat{X}_{s}|^2]\leq4e^{(8\gamma^2(T+\theta-t)+12m^2\gamma^2\rho^2)s}\leq4e^{8\gamma^2(T+\theta)^2+12m^2\gamma^2\rho^2(T+\theta)}.$$
Hence $E[\int^{T}_t|\hat{X}_{r}|^2dr]=\int^{T}_t E[|\hat{X}_{r}|^2]dr<+\infty$ and by Doob's martingale inequality, we deduce
\begin{align*}
&E[\sup_{s\in[t,T]}|\hat{X}_s|^2]\\
& \leq4+4E[\sup_{s\in[t,T]}|\int^s_ta_{r}\hat{X}_{r}dr|^2]+4E[\sup_{s\in[t,T]}|\int^s_t\mu_{r-\theta}\hat{X}_{r-\theta}dr|^2]\\
&+4E[\sup_{s\in[t,T]}|\int^s_t\hat{X}_{r-}b_{r-}
(\Psi_r^\dagger)'dM_{r}|^2]\\
& \leq4+8\gamma^2TE[\int^T_t|\hat{X}_{r}|^2 dr]+16E[|\int^T_t\hat{X}_{r-}b_{r-}
(\Psi_r^\dagger)'dM_{r}|^2].
\end{align*}
Similar to the above proof, we obtain $E[\sup_{s\in[t,T]}|\hat{X}_s|^2]<+\infty$.
\end{proof}
We return to the proof of Theorem \ref{bft1}. By Proposition \ref{prop} and Lemma \ref{Xbounded}, we know $E[\sup_{s\in[t,T']}|\tilde{L}|]<+\infty$. So
by Lemma \ref{local}, we deduce $\tilde{L}$ is a martingale. Because $\hat{X}_{t}=1$ and $\hat{X}_s=0$, $s\in [ t-\theta,t)$, taking conditional expectations under
$\mathcal{F}_{t},$ we have
$$\begin{array}{ll}
Y_{t}\\[0.2cm]
=E^{\mathcal{F}_{t}}[\hat{X}_TY_T+\int^T_{t} \hat{X}_s\mu_{s}E^{\mathcal{F}_{s}}[Y_{s+\theta
}]ds +\int^T_{t}\hat{X}_s\varphi_sds-\int^T_{t} Y_s\mu_{s-\theta}\hat{X}_{s-\theta}ds]\\[2mm]
=E^{\mathcal{F}_{t}}[\hat{X}_TY_T+\int^T_{t}\hat{X}_s\varphi_sds]+E^{\mathcal{F}_{t}}[\int^T_{t}(\hat{X}_s\mu_s Y_{s+\theta
}- \hat{X}_{s-\theta}\mu_{s-\theta}Y_s)ds]\\[0.2cm]
=E^{\mathcal{F}_{t}}[\hat{X}_TY_T+\int^T_{t}\hat{X}_s\varphi_sds]+E^{\mathcal{F}_{t}}[\int^{T+\theta}_{t+\theta}\hat{X}_{s-\theta}\mu_{s-\theta}Y_sds
-\int^{T}_{t}\hat{X}_{s-\theta}\mu_{s-\theta}Y_sds]\\[0.2cm]
=E^{\mathcal{F}_{t}}[\hat{X}_TY_T+\int^T_{t}\hat{X}_s\varphi_sds+\int^{T+\theta}_{T}\hat{X}_{s-\theta}\mu_{s-\theta} Y_{s
}ds
-\int^{t+\theta}_{t}\hat{X}_{s-\theta}\mu_{s-\theta} Y_{s
}ds]\\[0.2cm]
=E^{\mathcal{F}_{t}}[\hat{X}_TU_T+\int^T_{t}\hat{X}_s\varphi_sds+\int^{T+\theta}_{T}\mu_{s-\theta}\hat{X}_{s-\theta}U_sds],~~~\mbox{a.e.,~a.s.}
\end{array}$$
By Lemma 2.21 in Elliott \cite{elliott}, we obtain
$Y_t=E^{\mathcal{F}_{t}}[\hat{X}_TU_T+\int^T_{t}\hat{X}_s\varphi_sds+\int^{T+\theta}_{T}\mu_{s-\theta}\hat{X}_{s-\theta}U_s ds$, for any $t \in [0,T]$, a.s.$\mbox{}\hfill\Box$
\section{Comparison theorem of one-dimensional anticipated BSDEs with Markov chain model}
\indent The main idea of our proof comes from the proof of the comparison theorem for
anticipated BSDEs with Brownian motion noise in Peng and Yang \cite{PengYang}.\\
\indent Let $
(Y^{(1)}_{\cdot},Z^{(1)}_{\cdot}),~(Y^{(2)}_{\cdot},Z^{(2)}_{\cdot})$
be respectively the solutions of the following two one-dimensional
anticipated BSDEs:
$$
\left \{\begin{array}{ll}
-dY^{(j)}_t =f_j(t,Y^{(j)}_t,Z^{(j)}_t,Y^{(j)}_{t+\delta(t)})dt -Z^{(j)}_tdM_{t},&0 \leq t\leq T;\\[1mm]
 \hskip.575cmY^{(j)}_{t}=\xi^{(j)}_t,&T\leq t \leq T+K,
 \end{array}
\right.
$$
where $j=1,2.$
\begin{theorem}\label{bft2} Assume $\xi_{\cdot}^{(1)},\xi_{\cdot}^{(2)}\in L^2_{\mathcal{F}}(T,T+K;\mathbb{R})$, $\delta$ satisfies
$(i),(ii),$ and $f_1,f_2$ satisfy conditions such that the above two anticipated BSDEs have unique solutions. Suppose
 \begin{enumerate}
\item $f_1$ satisfies (H1), moreover, the Lipschitz constant $c_2$ of $f_1$ satisfies $$c_2\|\Psi_t^{\dagger}\|_{N \times N} \sqrt{6m}< 1, ~~~\text{ for any }~t \in [0,T],$$ where $\Psi$ is given in \eqref{Psi} and $m>0$ is the bound of $\|A_t\|_{N\times N}$, for any $t\in[0,T]$.
\item  for any $t\in[0,T],$ $y\in\mathbb{R},$ $z\in\mathbb{R}^N$,
$f_1(t,y,z,\cdot)$ is increasing, i.e., $f_1(t,y,z,\theta_r)$ $\geq
f_1(t,y,z,\theta'_r)$, if $\theta_r\geq\theta'_r$,
$\theta_{\cdot},\theta'_{\cdot}\in
L^2_{\mathcal{F}}(t,T+K;\mathbb{R}),r\in[t,T+K]$.\end{enumerate}
 If
$\xi^{(1)}_s\leq\xi^{(2)}_s,~s\in[T,T+K],$ and
$f_1(t,Y^{(2)}_t,Z^{(2)}_t,Y^{(2)}_{t+\delta(t)})\leq
f_2(t,Y^{(2)}_t,Z^{(2)}_t,Y^{(2)}_{t+\delta(t)}),$ a.e., a.s., then
$$P( Y_t^{(1)}\leq Y_t^{(2)},~~\text{ for any } t \in [0,T])=1.$$\end{theorem}
\noindent{\bf Proof.} Set
 $$ \left \{ \begin{array}{ll} Y^{(3)}_t = \xi_T ^{(1)} +\int ^T_tf_1(s,Y^{(3)}_s,Z^{(3)}_s,Y^{(2)}_{s+\delta(s)})ds  -\int ^T_t (Z^{(3)}_s)'d M_{s},&t\in[0,T]; \\[1mm]
 Y^{(3)}_{t} =
 \xi^{(1)}_t,&t\in[T,T+K].
\end{array} \right. $$
By Lemma \ref{existence}, we know there exists a solution $(Y^{(3)}_{\cdot},Z^{(3)}_{\cdot})\in L^2_{\mathcal{F}}(0,T;\mathbb{R})\times H^2(0,T;\mathbb{R}^N)$
to the above BSDE. Moreover, this solution is
unique up to indistinguishability for $Y.$ and equality $d\langle
M,M\rangle_t\times\mathbb{P}$-a.s. for $Z$. Set
$\tilde{f}_t=f_2(t,Y^{(2)}_t,Z^{(2)}_t,Y^{(2)}_{t+\delta(t)})-f_1(t,Y^{(2)}_t,Z^{(2)}_t,Y^{(2)}_{t+\delta(t)})$
and
$y_.=Y_.^{(2)}-Y_.^{(3)},~z_.=Z_.^{(2)}-Z_.^{(3)},~\tilde{\xi}_.=\xi_.^{(2)}-\xi_.^{(1)}.$
Then the pair $(y,z)$ can be regarded as the solution to the linear
BSDE
$$ \left \{ \begin{array}{ll}
y_t=\tilde{\xi}_T+\int ^T _t(a_sy_s+b_sz_s+\tilde{f}_s)ds-\int ^T
_tz_sdM_s,&t\in[0,T];\\[1mm]
y_t=\tilde{\xi}_t,&t\in[T,T+K],
\end{array} \right.
$$
where
$$\begin{array}{ll}
a_s=\begin{cases}\dfrac{f_1(t,Y^{(2)}_t,Z^{(2)}_t,Y^{(2)}_{t+\delta(t)})-f_1(t,Y^{(3)}_t,Z^{(2)}_t,Y^{(2)}_{t+\delta(t)})}{y_s},&\mbox{if } \,y_s\neq 0;\\
  0,& \mbox{if } \,y_s= 0,\end{cases} \\\\
b_s=\begin{cases}\dfrac{f_1(t,Y^{(3)}_t,Z^{(2)}_t,Y^{(2)}_{t+\delta(t)})-f_1(t,Y^{(3)}_t,Z^{(3)}_t,Y^{(2)}_{t+\delta(t)})}{|z_s|_N^2}z_s',& \mbox{if } \,z_s\neq 0;\\
 0,& \mbox{if } \,z_s= 0.\end{cases}
\end{array} $$
Since $f_1$ satisfies (H1), we deduce for any $s\in[0,T]$, $|a_s|\leq c_1$ and by Lemma \ref{normbound},
$$|b_s|_N\leq c_2\frac{\|z_s\|_{X_s}\cdot|z_s|_N}{|z_s|_N^2}\leq c_2\sqrt{3m}.$$
By Lemma \ref{Duality}, we know
 \begin{equation*}
P (y_t = E [ \tilde{\xi}_T U_T + \int_t^T \tilde{f}_s U_s ds | \mathcal{F}_t], ~ \text{for any}~ t\in [0,T])=1,
\end{equation*}
where $U$ is the solution of a one-dimensional SDE
\begin{equation}\label{sde}
\begin{cases}
dU_s= U_s a_s ds + U_{s-} b_{s-} (\Psi_s^{\dagger})'dM_s, ~~ s \in [t,T];\\
~~U_t =1.
\end{cases}
\end{equation}
Denote $$dV_s = a_sds + b_{s-} (\Psi_s^{\dagger})'dM_s,~~s\in[0,T].$$ The solution to SDE \eqref{sde} is given by the Dol\'{e}an-Dade exponential (See \cite{elliott}):
\[U_s = \exp(V_s - \frac{1}{2}\lel V^c,V^c\rir_s)\prod_{0 \leq u\leq s} (1+\Delta V_u )e^{-\Delta V_u},~~s\in[0,T],\]
 where $$\Delta V_u = b_{u-} (\Psi_u^{\dagger})' \Delta M_u = b_{u-} (\Psi_u^{\dagger})' \Delta X_u.$$
 Since $f_1$ satisfies $c_2\|\Psi_t^{\dagger}\|_{N \times N} \sqrt{6m}< 1,$ for any $t \in [0,T],$ where $\Psi$ is given in \eqref{Psi} and $m>0$ is the bound of $\|A_t\|_{N\times N}$, for any $t\in[0,T]$, by Lemma \ref{normbound} we have
\begin{align*}
|\Delta V_u| \leq |b_{u-} |_N \cdot \|(\Psi_u^{\dagger})'\|_{N\times N} \cdot |\Delta X_u|_N
< c_2\sqrt{3m}\dfrac{1}{\sqrt{6m}c_2}\sqrt{2}
 =1.
\end{align*}
Hence we have $U_s >0$, $s\in[0,T]$. As $\tilde{\xi}_T \geq 0$, a.s., and $\tilde{f}_s \geq 0$, a.e., a.s., we know for any $t \in [0,T]$,
 $$
y_t = E [ \tilde{\xi}_T U_T + \int_t^T \tilde{f}_s U_s ds |\mathcal{F}_t] \geq0, ~a.s.
$$
 Since $y_.$ is RCLL, by Lemma 2.21 in Elliott \cite{elliott}, we obtain
 \[P(Y^{(2)}_t\geq Y^{(3)}_t, ~ \text{for any } t \in [0,T])=P(y_t \geq 0, ~ \text{for any } t \in [0,T])=1.\]
Set
$$ \left \{\begin{array}{ll}
 Y^{(4)}_t = \xi_T ^{(1)} +\int^T_tf_1(s,Y^{(4)}_s,Z^{(4)}_s,Y^{(3)}_{s+\delta(s)})ds  -\int ^T_t
(Z^{(4)}_s)'dM_{s},&t\in[0,T];\\[1mm]
 Y^{(4)}_{t} =\xi^{(1)}_t,&t\in[T,T+K].
\end{array} \right. $$
Recall for any $t\in[0,T],y\in\mathbb{R},z\in\mathbb{R}^N,$
$f_1(t,y,z,\cdot)$ is increasing and $Y^{(2)}_t\geq Y^{(3)}_t,$ for any $ t \in [0,T]$, a.e. Also, $f_1$ satisfies $c_2\|\Psi_t^{\dagger}\|_{N \times N} \sqrt{6m}<1$ for $t \in [0,T]$. So
by Lemma \ref{CT'} we
obtain
$$P(Y^{(3)}_t\geq Y^{(4)}_t,~~\text{ for any } t \in [0,T])=1.$$
For $n=5,6,\cdot\cdot\cdot,$ we consider the following sequence of classical
BSDEs on Markov chain:
 $$ \left \{\begin{array}{ll}
 Y^{(n)}_t = \xi_T ^{(1)} +\int^T_tf_1(s,Y^{(n)}_s,Z^{(n)}_s,Y^{(n-1)}_{s+\delta(s)})ds  -\int ^T_t
(Z^{(n)}_s)'dM_{s},&t\in[0,T];\\[1mm]
 Y^{(n)}_{t} =\xi^{(1)}_t,&t\in[T,T+K].
\end{array} \right. $$
Similarly for any $n\in\mathbb{N},n\geq4$, we know the above equation has a unique solution $(Y^{(n)}_.,Z^{(n)}_. )\in L^2_{\mathcal{F}}(0,T;\mathbb{R})\times H^2(0,T;\mathbb{R}^N)$. Moreover, there exists a subset $A_n\subseteq\Omega$ with $P(A_n)=1$ such that for any $\omega\in A_n$, $Y_t^{(n)}(\omega)\geq Y_t^{(n+1)}(\omega),$ for any $t\in[0,T].$ Hence
$$
P(\bigcap_{n=4}^{+\infty}A_n) = 1 - P(\bigcup_{n=4}^{+\infty} A_n^c)
\geq 1 -\sum_{n=4}^{+\infty}P(A^c_n)
=1.
$$
That is,
$$P(Y_t^{(4)}\geq Y_t^{(5)}\geq \ldots\geq Y_t^{(n)}\geq\ldots, ~~\text{for
any }t\in[0,T])=1.$$
So
$$P(Y_t^{(2)}\geq Y_t^{(3)}\geq Y_t^{(4)}\geq Y_t^{(5)}\geq\ldots\geq Y_t^{(n)}\geq\ldots, ~~\text{for
any }t\in[0,T])=1.$$
Let $\beta>0$ be an arbitrary constant and $c=\max\{c_1,c_2\}$. We use
$\parallel\nu(\cdot)\parallel_{L^2}$ and $\parallel\mu(\cdot)\parallel_{\hat{H}^2}$ in the proof of Theorem
\ref{bft3000000} as the norms in the Banach spaces
$L^2_\mathcal{F}(0,T+K;\mathbb{R})$ and $\hat{H}^2(0,T+K;\mathbb{R}^{N})$, respectively. Set
$\hat{Y}^{(n)}_s=Y^{(n)}_s-Y^{(n-1)}_s,~\hat{Z}^{(n)}_s=Z^{(n)}_s-Z^{(n-1)}_s,~n\geq4.$
Then $(\hat{Y}^{(n)}_{\cdot},~\hat{Z}^{(n)}_{\cdot})$ satisfies the
following BSDE
 $$ \left \{
\begin{array}{ll} \hat{Y}^{(n)}_t = \int
^T_t(f_1(s,Y^{(n)}_s,Z^{(n)}_s,Y^{(n-1)}_{s+\delta(s)})-f_1(s,Y^{(n-1)}_s,Z^{(n-1)}_s,Y^{(n-2)}_{s+\delta(s)}))ds\\[0.1cm]
\hskip1.2cm-\int^T_t (\hat{Z}^{(n)}_s)'dM_{s},\hskip1cmt\in[0,T];\\[0.1cm]
 \hat{Y}^{(n)}_{t}=0,\hskip3.5cmt\in[T,T+K].
\end{array} \right. $$
Apply It\^{o}'s formula to $e^{\beta s} |\hat{Y}_s|^2$ for $s\in[0,T]$ and then take the expectation:
\begin{align*}
&E [  |\hat{Y}_0^{(n)}|^2 ]+ E [ \int_0^T \beta |\hat{Y}_s^{(n)}|^2 e^{\beta s} ds ]  + E [ \int_0^T\|\hat{Z}_s^{(n)}\|^2_{X_s}e^{\beta s} ds ] \\
& = 2 E [ \int_0^T \hat{Y}_s^{(n)} (f_1(s,Y^{(n)}_s,Z^{(n)}_s,Y^{(n-1)}_{s+\delta(s)})-f_1(s,Y^{(n-1)}_s,Z^{(n-1)}_s,Y^{(n-2)}_{s+\delta(s)})) e^{\beta s}ds ]\\
&\leq  E [ \int_0^T (\frac{\beta}{2}|\hat{Y}^{(n)}_s|^2+\frac{2}{\beta} |f_1(s,Y^{(n)}_s,Z^{(n)}_s,Y^{(n-1)}_{s+\delta(s)})-f_1(s,Y^{(n-1)}_s,Z^{(n-1)}_s,Y^{(n-2)}_{s+\delta(s)})|^2)e^{\beta s} ds ].
\end{align*}
Thus
\begin{align*}
&
E[\int^T_0(\dfrac{\beta}{2}|\hat{Y}^{(n)}_s|^{2}+\|\hat{Z}_s^{(n)}\|^2_{X_s})e^{\beta
s}ds] \\
&\leq
\dfrac{2}{\beta}E[\int^T_0|f_1(s,Y^{(n)}_s,Z^{(n)}_s,Y^{(n-1)}_{s+\delta(s)})-f_1(s,Y^{(n-1)}_s,Z^{(n-1)}_s,Y^{(n-2)}_{s+\delta(s)})|^{2}e^{\beta
s}ds]\\
&
\leq\dfrac{6c^2}{\beta}E[\int^T_0(|\hat{Y}^{(n)}_s|^2+\|\hat{Z}^{(n)}_s\|_{X_s}^2+|\hat{Y}^{(n-1)}_{s+\delta(s)}|^2)e^{\beta
s}ds]\\
&
\leq\dfrac{6c^2}{\beta}E[\int^T_0(|\hat{Y}^{(n)}_s|^2+\|\hat{Z}^{(n)}_s\|_{X_s}^2)e^{\beta
s}ds]+\dfrac{6c^2L}{\beta}E[\int^T_0|\hat{Y}^{(n-1)}_s|^2e^{\beta
s}ds].
\end{align*}
Set $\beta=18c^2L+18c^2+3$. Then
\begin{align*}\dfrac{2}{3}E[\int^T_0(|\hat{Y}^{(n)}_s|^{2}+\|\hat{Z}^{(n)}_s\|_{X_s}^2)e^{\beta
s}ds]&
 \leq\dfrac{1}{3}E[\int^T_0|\hat{Y}^{(n-1)}_s|^{2}e^{\beta s}ds]\\
&\leq\dfrac{1}{3}E[\int^T_0(|\hat{Y}^{(n-1)}_s|^{2}+\|\hat{Z}^{(n-1)}_s\|^{2}_{X_s})e^{\beta
s}ds].
\end{align*}
Hence
$$
E[\int^T_0(|\hat{Y}^{(n)}_s|^{2}+\|\hat{Z}^{(n)}_s\|_{X_s}^{2})e^{\beta
s}ds]
\leq(\frac{1}{2})^{n-4}E[\int^T_0(|\hat{Y}^{(4)}_s|^{2}+\|\hat{Z}^{(4)}_s\|^{2}_{X_s})e^{\beta
s}ds].
$$
It follows that $(Y^{(n)}_{\cdot})_{n\geq4}$ and
$(Z^{(n)}_{\cdot})_{n\geq4}$ are, respectively, Cauchy sequences in
$L^2_{\mathcal{F}}(0,T+K;\mathbb{R})$ and in
$\hat{H}^2(0,T+K;\mathbb{R}^{N})$. Denote their limits by
$Y_{\cdot}$ and $Z_{\cdot}$, respectively. Then $$P(Y_t^{(2)}\geq Y_t^{(3)}\geq \ldots\geq Y_t^{(n)}\geq\ldots\geq Y_t, ~~\text{for
any }t\in[0,T])=1.$$ Since
$L^2_{\mathcal{F}}(0,T+K;\mathbb{R})$ and $\hat{H}^2(0,T+K;\mathbb{R}^{N})$
are both Banach spaces, we obtain $(Y_{\cdot},Z_{\cdot})\in
L^2_{\mathcal{F}}(0,T+K;\mathbb{R})\times \hat{H}^2(0,T+K;\mathbb{R}^{N}).$
Note for any $\in[0,T]$,
$$\begin{array}{lll}
\hskip.44cmE[\int^T_t|f_1(s,Y^{(n)}_s,Z^{(n)}_s,Y^{(n-1)}_{s+\delta(s)})-f_1(s,Y_s,Z_s,Y_{s+\delta(s)})|^{2}e^{\beta
s}ds]\\[0.2cm]
\leq3c^2E[\int^T_t(|Y^{(n)}_s-Y_s|^2+\|Z^{(n)}_s-Z_s\|_{X_s}^2+L|Y^{(n-1)}_s-Y_s|^2)e^{\beta
s}ds]\rightarrow0,~n\rightarrow\infty.
\end{array}$$
Therefore, $(Y_{\cdot},Z_{\cdot})$ satisfies the following
anticipated BSDE
 $$ \left \{ \begin{array}{ll} Y_t = \xi_T ^{(1)} +\int ^T_tf_1(s,Y_s,Z_s,Y_{s+\delta(s)})ds  -\int ^T_t Z'_sdM_{s},&0 \leq t\leq T;\\[1mm]
 Y_{t} = \xi^{(1)}_t,&T\leq t \leq
 T+K.
\end{array} \right. $$
By Theorem \ref{bft3000000} we know
$$P(Y_t=Y^{(1)}_t, ~ \text{for any } t \in [0,T])=1.$$ Because $P(Y^{(2)}_t\geq
Y^{(3)}_t\geq Y_t,$ for any $t\in[0,T])=1$, it holds immediately that
$$\hskip3.5cm P(Y^{(1)}_t\leq
Y^{(2)}_t, ~ \text{for any } t \in [0,T])=1.\hskip3cm \mbox{}\hfill\Box$$

\end{document}